\documentclass[reqno,centertags, 12pt]{amsart}
\usepackage{amsmath,amsthm,amscd,amssymb}
\usepackage{latexsym}
\newcommand{\bbC}{{\mathbb{C}}}
\newcommand{\bbD}{{\mathbb{D}}}

\newcommand{\calC}{{\mathcal{C}}}

\newcommand{\calH}{{\mathcal H}}
\newcommand{\calK}{{\mathcal K}}
\newcommand{\calL}{{\mathcal L}}
\newcommand{\calM}{{\mathcal M}}

\newcommand{\bdone}{{\boldsymbol{1}}}
\newcommand{\dott}{\,\cdot\,}

\newcommand{\lb}{\label}
\newcommand{\f}{\frac}

\newcommand{\ol}{\overline}
\newcommand{\ti}{\tilde  }
\newcommand{\wti}{\widetilde  }

\newcommand{\Arg}{\text{\rm{Arg}}}

\newcommand{\ess}{\text{\rm{ess}}}

\newcommand{\supp}{\text{\rm{supp}}}

\newcommand{\bi}{\bibitem}

\newcommand{\beq}{\begin{equation}}
\newcommand{\eeq}{\end{equation}}
\newcommand{\ba}{\begin{align}}
\newcommand{\ea}{\end{align}}

\newcounter{smalllist}
\newenvironment{SL}{\begin{list}{{\rm\roman{smalllist})}}{%
\setlength{\topsep}{0mm}\setlength{\parsep}{0mm}\setlength{\itemsep}{0mm}%
\setlength{\labelwidth}{2em}\setlength{\leftmargin}{2em}\usecounter{smalllist}%
}}{\end{list}}

\DeclareMathOperator{\Real}{Re}

\allowdisplaybreaks
\numberwithin{equation}{section}
\newtheorem{theorem}{Theorem}[section]

\newtheorem*{p2.1}{Proposition 2.1}
\newtheorem{proposition}[theorem]{Proposition}
\newtheorem{lemma}[theorem]{Lemma}
\newtheorem{corollary}[theorem]{Corollary}
\theoremstyle{definition}

\newtheorem{example}[theorem]{Example}

\theoremstyle{remark}
\newtheorem*{remark}{Remark}
\newcommand{\abs}[1]{\lvert#1\rvert}

\begin{document}
\title[Rank One Perturbations and the Zeros of POPUC]
{Rank One Perturbations and the Zeros of
Paraorthogonal Polynomials on the Unit Circle}
\author[B. Simon]{Barry Simon*}

\thanks{$^*$ Mathematics 253-37, California Institute of Technology, Pasadena, CA 91125.
E-mail: bsimon@caltech.edu. Supported in part by NSF grant DMS-0140592}

\date{May 9, 2006}
%\keywords{orthogonal polynomials, Verblunsky coefficients, Szeg\H{o}'s theorem}
%\subjclass[2000]{42C05, 30E05, 42A70}

\begin{abstract} We prove several results about zeros of paraorthogonal polynomials
using the theory of rank one perturbations of unitary operators. In particular,
we obtain new details on the interlacing of zeros for successive POPUC.
\end{abstract}

\maketitle

%%%%%%%%%%%%%%%%%%%%%%%%%%%%%%%
\section{Introduction} \lb{s1}
%%%%%%%%%%%%%%%%%%%%%%%%%%%%%%%

This note concerns an aspect of the theory of orthogonal polynomials on
the unit circle (OPUC); for background, see \cite{Szb,GBk1,OPUC1,OPUC2,1Ft}.
Given a nontrivial probability measure, $d\mu$, on $\partial\bbD=\{z\in\bbC
\mid \abs{z}=1\}$, we let $\Phi_n(z)$ (we use $\Phi_n (z,d\mu)$ when $d\mu$
needs to be explicit) be the monic orthogonal polynomials. They obey the
Szeg\H{o} recursion relations
\begin{align}
\Phi_{n+1}(z) &= z\Phi_n(z) - \bar\alpha_n \Phi_n^*(z)  \lb{1.1} \\
\Phi_n^*(z) &= z^n \, \ol{\Phi_n (1/\bar z)} \lb{1.2}
\end{align}
where $\{\alpha_n\}_{n=0}^\infty\in\bbD^\infty$ are the Verblunsky coefficients.
$d\mu \leftrightarrow \{\alpha_n\}_{n=0}^\infty$ sets up a one-one correspondence
between $\bbD^\infty$ and nontrivial probability measures (Verblunsky's theorem).

Given $\beta\in\partial\bbD$, the paraorthogonal polynomials (POPUC) are
defined by ({\it Note\/}: \cite{OPUC1} uses $\beta$ where \eqref{1.3} uses
$-\bar\beta$; \eqref{1.3} is the right convention.)
\begin{equation} \lb{1.3}
\Phi_n(z,d\mu;\beta) = z\Phi_{n-1} (z,d\mu) - \bar\beta\Phi_{n-1}^*(z,d\mu)
\end{equation}
More generally, we will consider a sequence $\{\beta_n\}_{n=1}^\infty\in\partial\bbD$
and
\begin{equation} \lb{1.4}
\ti\Phi_n(z) = \Phi_n (z,d\mu;\beta_n)
\end{equation}

POPUC were introduced at least as early as Jones, Nj\r{a}stad, and Thron
\cite{JNT89}. About five years ago, Cantero, Moral, and Vel\'azquez
\cite{CMV02} and Golinskii \cite{Gol02} realized that zeros of POPUC shared
many properties of zeros of OPRL and independently proved a number of basic
results about these zeros. Cantero et al.\ \cite{CMV06} recently proved
additional results. The basic tool in \cite{CMV02,Gol02} is the Christoffel--Darboux
formula; \cite{CMV06} also exploits the CMV matrix. Our goal in this paper is to use
the theory of rank one perturbations of unitary matrices to recover many of the
basic results about zeros of POPUC and prove some new results. In particular, we will
illuminate the issue of interlacing of the zeros of successive POPUC.

First, some notation. Given distinct $z,w\in\partial\bbD$, $(z,w)$ is the set of
points, $\zeta$, in $\partial\bbD$ with
\begin{equation} \lb{1.5}
\Arg (z) < \Arg(\zeta) < \Arg (w)
\end{equation}
where a branch of $\Arg$ is chosen so $0<\Arg (w)-\Arg(z) < 2\pi$. An ordered
set of points $(z_1, \dots, z_\ell)\in\partial\bbD^\ell$ is called {\it cyclicly
ordered} if each $(z_j,z_{j+1})_{j=1}^\ell$ and $(z_\ell, z_1)$ contain no other
$z_j$'s. The ordering is fixed by such cyclicity up to a single choice. We will always
assume zeros of POPUC are cyclicly ordered.

Two cyclicly ordered sets $(z_1, \dots, z_\ell)$ and $(w_1, \dots, w_\ell)$ in
$\partial\bbD^\ell$ are said to {\it strictly interlace} if after a cyclic
permutation of the $w$'s, we have $w_j\in (z_j, z_{j+1})$, $j=1,2,\dots, \ell-1$,
$w_\ell\in (z_\ell, z_1)$. This, of course, implies $z_j\in (w_{j-1}, w_j)$, $j=2,3,
\dots, \ell$ and $z_1\in (w_\ell, w_1)$.

For $\{\alpha_j\}_{j=0}^\infty$, the second kind polynomials, $\Psi_n (z,d\mu)$ are
defined, as usual, to be the $\Phi_n$'s associated to $\ti\alpha_j = -\alpha_j (d\mu)$.
We define
\begin{equation} \lb{1.6}
\Psi_n (z,d\mu;\beta) =z\Psi_{n-1} (z,d\mu) -\bar\beta \Psi_{n-1}^*(z,d\mu)
\end{equation}

We can now state our main results:

\begin{theorem}[\cite{CMV02,Gol02}]\lb{T1.1} If $(w_0, w_1)$ is an interval disjoint
from $\supp(d\mu)$, then for any choice of $\beta$ and any $n$, $\Phi_n (z,d\mu;\beta)$
has at most one zero in $(w_0, w_1)$.
\end{theorem}

The following has also been proven by Wong \cite{Wong}:

\begin{theorem}\lb{T1.2}  Let $(z_1, \dots, z_n)$ be the zeros of some $\Phi_n (z,d\mu;\beta)$
and $(w_1, \dots, w_n)$ of $\Psi_n (z,d\mu; -\beta)$. {\rm{(}}Note: $-\beta$, not
$\beta$.{\rm{)}} Then these zeros strictly interlace.
\end{theorem}

\begin{theorem}[\cite{CMV02,Gol02}]\lb{T1.3} Fix $d\mu$ and $n$ and distinct $\beta, \beta'$
in $\partial\bbD$. Then the zeros of $\Phi_n (z,d\mu;\beta)$ and $\Phi_n (z,d\mu;\beta')$
strictly interlace.
\end{theorem}

The power of our approach is shown by the refined version we obtain relating zeros
of $\ti\Phi_{n+1}$ and $\ti\Phi_n$. We will need the following computed sequence
in $\partial\bbD$:
\begin{equation} \lb{1.7}
\lambda_n = \bar\beta_{n+1} \bar\beta_n \biggl( \f{\beta_n\alpha_n -1}
{\bar\beta_n \bar\alpha_n -1}\biggr)
\end{equation}

\begin{theorem}\lb{T1.4} For each $n$, one of two possibilities holds:
\begin{SL}
\item[{\rm{(i)}}] $\ti\Phi_n$ and $\ti\Phi_{n+1}$ have no zeros in common. In that case,
$\lambda_n$ is not a zero of either, and $\{$zeros of $\ti\Phi_n\}\cup\{\lambda_n\}$
strictly interlace $\{$zeros of $\ti\Phi_{n+1}\}$.

\item[{\rm{(ii)}}] $\ti\Phi_n$ and $\ti\Phi_{n+1}$ have a single zero in common.
In that case, $\lambda_n$ is that zero and $\{$zeros of $\ti\Phi_n\}$ strictly
interlace $\{$zeros of $\ti\Phi_{n+1}\}\backslash\{\lambda_n\}$.
\end{SL}
\end{theorem}

\begin{corollary}\lb{C1.5} If $\ti\Phi_1, \ti\Phi_2,\ti\Phi_3, \dots$ have a common
zero at $\lambda$, then $\beta_n$ are given inductively by
\begin{align}
\beta_1 &= \bar\lambda \lb{1.8} \\
\beta_{n+1} &= \bar\lambda \bar\beta_n \biggl( \f{\beta_n \alpha_n-1}
{\bar\beta_n \bar\alpha_n -1}\biggr)
\end{align}
\end{corollary}

\begin{example}\lb{E1.6} $\alpha\equiv 0$. Then $\beta_n =\bar\lambda^n$ and $\ti\Phi_n(z)
=z^n -\lambda^n$ precisely the POPs with a zero at $\lambda$ for all $n$.
\qed
\end{example}

The key to our proofs is the connection of $\ti\Phi_n$ to CMV matrices
{\cite{CMV, OPUC1, Evansproc}. $\ti\Phi_n$ is the determinant of a suitable
finite CMV matrix, and so its zeros are the eigenvalues. All our results concern
what happens to eigenvalues of unitary matrices under rank one perturbations.
Section~\ref{s2} discusses general rank one perturbations, and Section~\ref{s3}
the application to POPUC.

\smallskip
It is a pleasure to thank Mar\'a Jos\'e Cantero and Lilian Wong for
useful discussions.

%%%%%%%%%%%%%%%%%%%%%%%%%%%%%%%%%%%%%%%%
\section{Rank One Perturbations} \lb{s2}
%%%%%%%%%%%%%%%%%%%%%%%%%%%%%%%%%%%%%%%%

Rank one perturbations of unitaries are discussed in Sections~1.3.9,
1.4.16, 3.2, 4.5, and 10.16 of \cite{OPUC1,OPUC2}, and some of the
results in this section are spread through that material.

If $U$ and $V$ are two unitaries on a finite- or infinite-dimensional
Hilbert space and $U-V$ is rank one, we pick a unit vector $\varphi\in
\ker (U-V)^\perp$ and note there must be a $\lambda\in\partial\bbD$ with
\begin{equation} \lb{2.1}
V\varphi =\lambda U\varphi
\end{equation}
and thus
\begin{equation} \lb{2.2}
V-U =(\lambda -1) \langle\varphi, \dott\rangle U\varphi
\end{equation}

It is convenient to define for $z\in\bbD$ and $A$ unitary:
\begin{align}
F_{A,\varphi}(z) &= \biggl\langle \varphi, \,\f{A+z}{A-z}\, \varphi\biggr\rangle \lb{2.3} \\
f_{A,\varphi}(z) &= z^{-1} (1-F(z))(1+F(z))^{-1} \lb{2.4}
\end{align}
$F$ is a Carath\'eodory function ($\Real\, F(z)>0$ on $\bbD$, $F(0)=1$) and $f$ a
Schur function ($\abs{f(z)}<1$ on $\bbD$). The spectral measure for $A,\varphi$ is given by
\begin{equation} \lb{2.5}
F_{A,\varphi} (z) = \int \f{e^{i\theta}+z}{e^{i\theta}-z}\, d\mu_{A,\varphi}(z)
\end{equation}
It is not hard to see that

\begin{proposition}\lb{P2.1} Let $\varphi$ be a cyclic vector for a unitary $A$
{\rm{(}}i.e., $\{A^k\varphi\}_{k=-\infty}^\infty$ is total{\rm{)}}. An interval $(w_0,w_1)$
in $\partial\bbD$ is disjoint from $\sigma_\ess (A)$ if and only if $f$ has an
analytic continuation through $(w_0,w_1)$ with $\abs{f(z)}=1$ on that interval.
In that case,
\begin{SL}
\item[{\rm{(a)}}] $\Arg(f)$ is strictly monotone increasing on $(w_0,w_1)$.
\item[{\rm{(b)}}] The only spectra of $A$ on $(w_0,w_1)$ are simple eigenvalues
precisely at the points $z$ where
\begin{equation} \lb{2.6}
zf(z)=1
\end{equation}
\end{SL}
\end{proposition}

$\Arg(f)$ is increasing since $\abs{f(z)}<1$ in $\bbD$ and $\abs{f(z)}=1$ on
$(w_0, w_1)$ implies $\partial\abs{f(re^{i\theta})}/\partial r \geq 0$ on
$(w_0,w_1)$. So by the Cauchy--Riemann equations, $\partial \Arg(f(e^{i\theta}))/
\partial\theta\geq 0$. (b) holds since
\begin{equation} \lb{2.7}
F(z) = \f{1+zf(z)}{1-zf(z)}
\end{equation}
has poles at points where \eqref{2.6} holds.

When \eqref{2.2} holds, a direct calculation (see (1.4.90) and the end of
Section~3.2 in \cite{OPUC1}) shows that

\begin{proposition}\lb{P2.2} If \eqref{2.2} holds, then
\begin{equation} \lb{2.8}
f_{V,\varphi}(z) = \lambda^{-1} f_{U,\varphi}(z)
\end{equation}
\end{proposition}

We immediately have

\begin{theorem}\lb{T2.3} Let \eqref{2.2} hold. If $(w_0,w_1)\cap\sigma(U)=\emptyset$,
then $V$\! has at most one eigenvalue in $[w_0,w_1]$ and no other spectrum there.
\end{theorem}

\begin{proof} Let $\calK$ be the cyclic subspace for $U$ and $\varphi$. Since $U=V$
on $\calK^\perp$ which is invariant for both, we can suppose $\varphi$ is cyclic.
In that case, picking $z_0\in (w_0, w_1)$ and then $\Arg(z(f(z)))$ so
$\Arg(z_0 f_0(z_0))\in (0,2\pi)$, we see $\Arg (z(f(z))\in (0,2\pi)$ on
all of $(w_0,w_1)$ since \eqref{2.6} has no solution there. By the strict
monotonicity of $\Arg(f)$, $zf(z)=\lambda$ has at most one solution in
$[w_0,w_1]$, so by Propositions~\ref{P2.1} and \ref{P2.2}, $V$ that
at most one eigenvalue there.
\end{proof}

\begin{proposition}\lb{P2.4} Let $U,V$ be unitaries on $\bbC^n$ so \eqref{2.2}
holds for $\lambda\neq 1$ and for $\varphi$ cyclic for $U$\!. Then the eigenvalues
of $U$\! and $V$\! strictly interlace.
\end{proposition}

\begin{proof} Since $U$ has a cyclic vector, its spectrum is simple so $zf(z)=1$
has $n$ solutions. Since $\Arg(z(f))$ is strictly monotone, $zf(z)=\lambda$ has $n$
solutions which interlace the solutions of $zf(z)=1$.
\end{proof}

One can say something about the case where $\varphi$ is not cyclic.

\begin{proposition}\lb{P2.5} Let $U,V$\! be unitaries on $\bbC^n$ so \eqref{2.2}
holds. Let $z_0,z_1$ be two eigenvalues of $U$\!. Then $V$\! has an eigenvalue in
$[z_0,z_1]$ {\rm{(}}$=(z_0, z_1)\cup\{z_0,z_1\}${\rm{)}}.
\end{proposition}

\begin{proof} Let $\calK$ be the cyclic subspace of $(U,\varphi)$ which is
invariant for $U$. If $z_0$ and $z_1$ are eigenvalues of $U\restriction\calK$,
$V$ has an eigenvalue in $(z_0,z_1)$ by Proposition~\ref{P2.4}. If not, since
$U\restriction \calK^\perp =V\restriction\calK^\perp$, either $z_0$ or $z_1$
is an eigenvalue of $V$\!.
\end{proof}

Finally, we have a specialized result that is precisely what we need to
prove Theorem~\ref{T1.4}:

\begin{proposition}\lb{P2.6} Let $U=U_1\oplus U_2$ on $\calK_1 \oplus \calK_2$,
two finite-dimensional subspaces of $\calH$, a space of dimension $n$. Let
$\varphi_j$ {\rm{(}}$j=1,2${\rm{)}} be cyclic vectors for $U_j$ on $K_j$.
Let $\varphi = a\varphi_1 \oplus b\varphi_2$ where $(a,b)\neq (0,0)$ and
$\abs{a}^2 + \abs{b}^2 =1$. Let $V$ be given by \eqref{2.2} with $\lambda \neq 1$.
If $U_1$ and $U_2$ have $\ell$ eigenvalues in common, then $V$\! has these $\ell$
common values as eigenvalues and its other $n-\ell$ eigenvalues strictly interlace
those of $U$.
\end{proposition}

\begin{proof} Since $\varphi_j$ is cyclic for $U_j$, any simple eigenvalue of
$U$ is in the cyclic subspace generated by $U,\varphi$. Moreover, any common
eigenvalue is a simple eigenvalue for $U\restriction\calK$ where $\calK=$
cyclic subspace of $\varphi$. Thus $U\restriction\calK$ has all the eigenvalues
of $U$ but with multiplicity $1$. The eigenvalues of $V\restriction\calK$
strictly interlace by Proposition~\ref{P2.4}. The eigenvalues of $V\restriction
\calK^\perp =U\restriction\calK^\perp$ are exactly the common eigenvalues.
\end{proof}

\begin{remark} $\varphi$ is cyclic if and only if $\ell=0$.
\end{remark}

%%%%%%%%%%%%%%%%%%%%%%%%%%%%%%%%%%%%%%%%%%%%%%%%%%%%%%%%
\section{Zeros of POPUC and Finite CMV Matrices} \lb{s3}
%%%%%%%%%%%%%%%%%%%%%%%%%%%%%%%%%%%%%%%%%%%%%%%%%%%%%%%%

Given a sequence $\{\gamma_n\}_{n=0}^\infty$ of elements in $\ol{\bbD}$,
one defines the CMV matrix $\calC (\{\gamma_n\}_{n=0}^\infty)$ on $\ell^2$ by
\begin{align}
\calC &= \calL\calM \lb{3.1} \\
\calL &= \Theta(\gamma_0) \oplus\Theta(\gamma_2)\oplus\cdots \lb{3.2} \\
\calM &= \bdone_{1\times 1} \oplus\Theta(\gamma_1)\oplus\Theta(\gamma_3)\oplus\cdots \lb{3.3}
\end{align}
where $\bdone_{1\times 1}$ is the one-dimensional identity matrix, and $\Theta$ is given by
\begin{align}
\Theta(\gamma) &= \left( \begin{array}{rr}
\bar\gamma & \tau \\
\tau & -\gamma
\end{array} \right) \lb{3.4} \\
\tau &= (1-\abs{\gamma}^2)^{1/2} \lb{3.5}
\end{align}

It is a fundamental result of Cantero, Moral, and Vel\'azquez \cite{CMV}
(see also \cite[Section~4.2]{OPUC1}) and see \cite{Evansproc} for other
references) that if $d\mu$ is a nontrivial probability measure on $\partial\bbD$,
$\chi_n$ is the basis of $L^2 (\partial\bbD,d\mu)$ obtained by applying
Gram--Schmidt to $1,z,z^{-1}, z^2, z^{-2}, \dots$, and $\alpha_n (d\mu)$ are
the Verblunsky coefficients of $d\mu$, then $\calC(\{\alpha_n (d\mu)\}_{n=0}^\infty)$
is the matrix of multiplication by $z$ on $L^2 (\partial\bbD,d\mu)$ in $\chi_n$ basis.
Note in this case that $\gamma_n\in\bbD$ (rather than some $\gamma_n\in\partial\bbD$),
in which case we call $\calC$ a proper CMV matrix.

If $\abs{\gamma}=1$, then $\tau =0$, and $\Theta(\gamma)$ is a direct sum of two
$1\times 1$ matrices, and so, if $\abs{\gamma_{n-1}}=1$, $\calC$ breaks into a direct sum
of an $n\times n$ matrix and an infinite piece. The finite piece, $\calC_n
(\{\gamma_0, \dots, \gamma_{n-1}\})$, is called a finite CMV matrix. It is not
hard to show that (see, e.g., \cite{Evansproc}):

\begin{proposition}\lb{P3.1} If $\gamma_n\in\bbD$ for all $n$, then $\delta_0\equiv
(1,0,\dots)^t$ is a cyclic vector for $\calC (\{\gamma_n\}_{n=0}^\infty)$.
If $\gamma_0, \dots, \gamma_{n-2}\in\bbD$, $\gamma_{n-1}\in\partial\bbD$,
then $\delta_0$ is a cyclic vector for $\calC_n (\{\gamma_m\}_{m=0}^{n-1})$.
\end{proposition}

Moreover (see \cite[Section~4.2]{OPUC1}),

\begin{proposition}\lb{P3.2} If $\alpha_0, \dots, \alpha_{n-2}\in\bbD$ and
$\beta\equiv \alpha_{n-1}\in\partial\bbD$, then
\begin{equation} \lb{3.6}
\Phi_n (z,d\mu_\alpha;\beta)=\det (z-\calC_n (\{\alpha_j\}_{j=0}^{n-1}))
\end{equation}
In particular, the zeros of $\ti\Phi_n$ are the eigenvalues of a finite CMV matrix.
\end{proposition}

Finally, we need the following, which generalizes Lemma~4.5.1 of \cite{OPUC1}:

\begin{lemma}\lb{L3.3} Let $\alpha\in\bbD$ and $\beta\in\partial\bbD$. Then
\begin{equation} \lb{3.7}
\Theta(\alpha) - \begin{pmatrix}
\beta & 0 \\ 0 & x \end{pmatrix}
\end{equation}
is rank one if and only if
\begin{equation} \lb{3.8}
x = \bar\beta \biggl( \f{\beta\alpha -1}{\bar\beta\bar\alpha -1}\biggr)
\end{equation}
\end{lemma}

\begin{proof} A $2\times 2$ matrix is rank one if and only if $\det(A)=0$. Since
\begin{equation} \lb{3.9x}
\det\biggl( \Theta(\alpha) - \begin{pmatrix} \beta & 0 \\ 0 & x \end{pmatrix}\biggr)
= (\bar\alpha - \beta)(-\alpha -x) -(1-\abs{\alpha}^2)
\end{equation}
we see \eqref{3.7} is rank one if and only if $\text{RHS of \eqref{3.9x}}=0$, which is
solved by \eqref{3.8}.
\end{proof}

{\it Note\/}: $\abs{x}=1$, so $\left(\begin{smallmatrix} \beta & 0 \\ 0 & x\end{smallmatrix}
\right)$ is unitary.

\begin{proof}[Proof of Theorem~\ref{T1.1}] Let $\calC$ be the CMV matrix, $\calC
(\{\alpha_n (d\mu)\}_{n=0}^\infty)$, of $d\mu$. Given $n$ and $\beta$, pick
$x$ so $\Theta (\alpha_{n-1})-\left(\begin{smallmatrix} \beta & 0 \\ 0 & x
\end{smallmatrix}\right)$ is rank one, and let $\ti\calC$ be the matrix obtained from
$\calC$ by replacing $\Theta (\alpha_{n-1})$ by $\left(\begin{smallmatrix}
\beta & 0 \\ 0 & x\end{smallmatrix}\right)$. Then $\ti\calC$ is unitary (by the
note after Lemma~\ref{L3.3}) and $\calC-\ti\calC$ is rank one. Thus, by
Theorem~\ref{T2.3}, $\ti\calC$ has at most one eigenvalue in $(w_0,w_1)$. But
$\ti\calC$ is a direct sum of $\calC_n (\{\alpha_0, \dots, \alpha_{n-2},\beta\})$
and another matrix, so $\calC_n$ has at most one eigenvalue in $(w_0,w_1)$. By
Proposition~\ref{P3.2}, zeros of $\ti\Phi_n$ are eigenvalues of $\calC_n$.
\end{proof}

\begin{proof}[Proof of Theorem~\ref{T1.2}] Let $\ti\alpha_{n-1}\equiv\beta$. By
Theorem~5.2 of \cite{Evansproc} (see also \cite[Theorem~4.2.9]{OPUC1}), $\calC_n
(\{-\ti\alpha_m\}_{m=0}^{n-1})$ is unitarily equivalent to $\ti\calC_n\equiv
\calL (\{\ti\alpha_m\}_{m=0}^{n-1})\wti\calM (\{\ti\alpha_m\}_{m=0}^{n-1})$
where $\wti\calM$ differs from $\calM$ by having $-\bdone_{1\times 1}$ in place
of $\bdone_{1\times 1}$. Thus, $\calC_n (\{\alpha_m\}_{m=0}^{n-1})-\ti\calC$ is
rank one, and Theorem~\ref{T1.2} follows from Propositions~\ref{P2.4} and
\ref{P3.1}.
\end{proof}

\begin{proof}[Proof of Theorem~\ref{T1.3}] If
\begin{gather}
\alpha_j =\alpha'_j =\alpha_j (d\mu) \qquad j=0,\dots, n-2  \lb{3.9} \\
\alpha_{n-1}=\beta \qquad \alpha'_{n-1}=\beta' \lb{3.10}
\end{gather}
then $\calC(\{\alpha_m\}_{m=0}^{n-1})-\calC(\{\alpha'_m\}_{m=0}^{n-1})$ is obviously
rank one. Moreover, $\delta_{n-1}$ is a cyclic vector since $\calC_n$ run backwards
is essentially another $\calC_n$ (with the initial $\bdone_{1\times 1}$ replaced by
$\bar\alpha_{n-1} \bdone_{1\times 1}$) or $\calC_n^t$. Thus, Theorem~\ref{T1.3} follows
from Proposition~\ref{P2.4}.
\end{proof}

\begin{proof}[Proof of Theorem~\ref{T1.4}] Let $\calC_{n+1}$ (resp.\ $\calC_n)$ be the
$(n+1)\times (n+1)$ (resp.\ $n\times n$) finite CMV matrix whose characteristic polynomial
is $\ti\Phi_{n+1}$ (resp.\ $\ti\Phi_n$). By Lemma~\ref{L3.3}, a rank one perturbation
turns $\calC_{n+1}$ into $\calC_n \oplus\lambda_n\bdone_{1\times 1}$ where $\lambda_n$
is given by \eqref{1.7}. The vector in the perturbation is $a\delta_{n-1} + b\delta_n$,
so Proposition~\ref{P2.6} applies and proves Theorem~\ref{T1.4}.
\end{proof}

As a final result:

\begin{theorem}\lb{T3.4} Let $m>n$. Then strictly between any pair of zeros of
$\ti\Phi_n$ is a zero of $\ti\Phi_m$.
\end{theorem}

\begin{proof} Let $\calC_n,\calC_m$ be as in the last proof. By a rank one perturbation,
$\calC_m$ can be changed to $\calC_n \oplus Q_{m-n}$. Now apply Proposition~\ref{P2.6}.
\end{proof}

\bigskip
%%%%%%%%%%%%%%%%%%%%%%%%%%%%%

\end{document}